\documentclass[12pt]{article}
\usepackage{geometry}                
\usepackage{graphicx}
\usepackage{amssymb}

\newtheorem{Theorem}{Theorem}[section]
\newtheorem{Proposition}[Theorem]{Proposition}

\newtheorem{Remark}[Theorem]{Remark}
\newtheorem{Definition}[Theorem]{Definition}

\newcommand{\bi}{\begin{enumerate}}
\newcommand{\ei}{\end{enumerate}}
\newcommand{\be}{\begin{equation}}
\newcommand{\ee}{\end{equation}}
\newcommand{\ba}{\begin{array}}
\newcommand{\ea}{\end{array}}

\def\ra{\rightarrow}

\def\t{\mathbf{t}}
\def\c{\mathbb{C}}

\def\q{\mathbb{Q}}
\def\z{\mathbb{Z}}

\def\cp{{\mathcal P}}

\title{First Chern class and birational germs of Kato surfaces}
\author{Massimiliano Pontecorvo}
\date{\today}

\usepackage{pst-func}    
\begin{document}
\maketitle

\begin{abstract}
We describe some relations between coefficients of irreducible components of the first Chern class \cite{fp15} and birational germs introduced by Dloussky \cite{dl16} for intermediate Kato surfaces.
\end{abstract}

\section{Introduction}

We study in this work some recent new aspects in the geometry of Kato surfaces after the works of Dloussky \cite{dl16} and Fujiki-Pontecorvo \cite{fp15}. In our terminology \it surface \rm precisely means a compact connected complex manifold of (complex) dimension $2$ and we will be concerned with the non-K\"ahler case which in this dimension is equivalent to ask for the first Betti number $b_1$ to be odd.

We further specialize to surfaces of non-K\"ahler type and Kodaira dimension $-\infty$ which form the class VII in Kodaira classification and are known to always satisfy $b_1=1$ \cite{bpv}. The subclass of such \it minimal \rm  surfaces (i.e. free from smooth rational curves of self-intersection $-1$) is denoted by VII$_0$  and finally we will be concerned with surfaces of class VII$_0^+$ meaning that the second Betti number $b_2$ is assumed to be positive. The classification of class-VII surfaces is complete only for  $b_2=0$ which consists of Hopf and Inoue-Bombieri surfaces; on the other hand the only known examples in class-VII$_0^+$ are surfaces with a global spherical shell (GSS) constructed by Ma. Kato \cite{ka77} for every $b_2\geq 1$.

We will concentrate on some geometric aspects of a subclass of Kato surfaces called \it intermediate, \rm let us explain why.  Every Kato surface $S$ has algebraic dimension zero and therefore a finite number of (compact complex) curves, in fact it always admits exactly $b_2(S)$-rational curves which we will denote $D_i$ for $i=1,\cdots, b:=b_2(S)$.

Kato's construction can be summarized in the following way which displays both the rational curves and the global spherical shell: one can start by blowing up the origin of a ball $0\in B \subset \c^2$ at $b$ points which are infinitely near and obtain in this way a non-compact surface $\tilde{B}$ with $b$ rational curves of negative self-intersection the last one of which $C_b$ being the only one of self-intersection $-1$. One then obtains a compact \it minimal \rm surface by considering a biholomorphism $\psi$ between the closure $\bar B \subset \c^2$ and a compact neighborhood of $p\in C_b\subset \tilde{B}$; the final step consists in identifying the two connected components of the boundary of $\tilde{B}\setminus \psi(B)$.

Topologically, $S$ is obtained from $\tilde{B}$ by a self-connected sum and therefore contains a global spherical shell (GSS) -- i.e. an open set $U$ biholomorphic to a neighborhood of $S^3\subset\c^2$ such that $S\setminus U$ is connected -- and the fundamental group satisfies $\pi_1(S)=\z$.
As a complex surface $S$ is always minimal because $C_b$ was the only curve in $\tilde{B}$ of self-intersection $-1$; more precisely we would like to notice at this point that most of the rational curves in $S$ will have self-intersection $-2$. Indeed, in order to have a curve $C$ with $C^2=-(k+2)\leq -3$ one needs to repeatedly blow up a fixed intersection point of previously created exceptional curves forcing $C$ to come along with a chain of $-2$-curves of lenght $(k-1)$. 

The following picture is meant to display this duality between rational curves of self-intersection $-(k+2)\leq -3$ which we indicate by a black node and chains of $(-2)$-curves of length $(k-1)$ indicated by a sequence of white nodes.

\medskip                
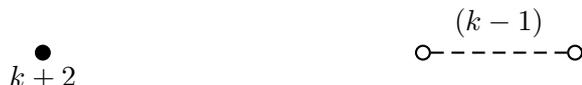
\begin{pspicture}(-6,-1)(6,1)


\cnode(2,0){3 pt}{A}
\cnode(4,0){3 pt}{B}


\ncline[linestyle=dashed]{-}{A}{B}
\Aput{\small{$(k-1)$}}


\qdisk(-3,0){3 pt}
\uput[270](-3,0){\small{$k+2$}}

\end{pspicture}
  \medskip

If $D=D_1+\dots +D_b$  is the reduced divisor of all rational curves on $S$, we will consider its \it dual graph \rm $\Gamma_D$ representing the configuration of these curves in $S$. $\Gamma_D$ consists of $b$ nodes each of which represents a curve with edges connecting those curves which do intersect each other. As it turns out this graph is always fairly simple and in the next section we will show how it is constructed; before we get to the subject of our work it seems appropriate to recall a finer subdivision of Kato surfaces in terms of (the dual graph of) $D$.
\begin{enumerate}
\item \it Enoki surfaces. \rm $D$ consists of a cycle of $(-2)$-curves, in other words $D^2=0$. $\Gamma_D$ is then a cycle of $b$ white nodes and $b$ edges.
\item \it Hyperbolic and half-Inoue surfaces. \rm These surfaces are also called Inoue-Hirzebruch surfaces, $D$ consists of one or two cycles with self-intersection $D^2=-b<0$; therefore black nodes appear in $\Gamma_D$.
\item \it Intermediate Kato surfaces. \rm In this case the dual graph $\Gamma_D$ is always connected and consists of a unique (non-empty) cycle with a positive number of trees appended to it.
 \end{enumerate}
 
 The name \it intermediate \rm is due to the Dloussky number dl$(S)$ of Kato surfaces which is defined to be the sum of the opposite intersection numbers of its rational curves so that the minimum is attained by Enoki surfaces for which dl$(Enoki)=2b$; in general, the duality between black nodes and chains of $(-2)$-curves easily implies that dl$(S)\leq 3b$ with equality holding only when there are no trees in $\Gamma_D$. Therefore $2<dl(S)<3$ applies precisely for the so called \it intermedate \rm case.
 
 The aim of the present work is to study intermediate Kato surfaces and try to relate some recent work of Dloussky \cite{dl16} on new birational structures with the results of \cite{fp15} in which we described the anti-canonical class $c_1=-K\in H^2(S,\z)$ as a linear combination of the irreducible components of $D$.
 
 In section 2 we will start by displaying the dual graph of such an intermediate surface, recall some properties of  multilinear forms which we defined in \cite{fp15} and show how they fit in the work of Dloussky in order to identify some topological invariants of $D\subset S$.
 
 In section 3 our interest moves to new birational germs found by Dloussky whose structure reflects the configuration of the irreducible components of the maximal curve $D\subset S$; we will show how the description of the first Chern class $c_1(S)\in H^2(S,\z)$ \cite{fp15} is related to some invariants of the germ of $S$ and in the last section we will study deformations of the pair $(S,D)$ which (by construction) are always parametrized by Dloussky's birational germ.

\section{Intermediate Kato surfaces.}

The following procedure to construct the dual graph $\Gamma_D$ of a Kato surface $S$ with reduced divisor of all rational curves  $D=D_1+\cdots +D_b$ was introduced by Dloussky \cite{dl84}. One may start with a \it simple \rm  Dloussky sequence as described in \cite[p.335]{ot08} or \cite[p.94]{dl16} of the following form
\be   \mathrm{Dl}S=[s_{k_1}s_{k_2}...s_{k_N}r_l]  \label{simple} \ee
It is a sequence of $b_2(S)=:b=k_1+k_2+\cdots+k_N+l$ positive integers each of which equals $-D_i^2$ for $i=1,\dots,b$; we also assume existence of at least one \it singular sequence \rm $s_k=(k+2,2,...,2)$ and exactly one \it regular sequence \rm $r_l=(2,...,2)$ and recall that each positive integers $k_j$ and $l$ represents the length of the corresponding subsequence. On the other hand it will be clear from the prescription below that a Dloussky sequence with $l=0$ corresponds to a Inoue-Hirzebruch surface while those with all $k_j=0$ represent Enoki surfaces.
  
The simple Dloussky sequence above is equivalent to say that the Kato surface $S$ is of intermediate type with only one tree and this can be seen as follows. 
Given the sequence $ \mathrm{Dl}S$ its dual graph is constructed by connecting an entry with value $a$ to the entry following $a-1$ places after it on the right, in cyclic order. There are exactly $b_2(S)$ rational curves in $D$ each of them can be assumed to be smooth (after taking a double covering, if necessary); correspondingly there are $b=b_2(S)$ nodes in $\Gamma_D$.  As already noticed, most curves have self-intersection number $-2$ and we will indicate them by a white node without any further reference to their self-intersection number. Furthermore, there are exactly $N$ curves whose self-intersection number is not greater than $-3$; we denote them by a black node and a positive integer indicating their opposite self-intersection number.

When $N$ is even, the black nodes of $S$ are evenly distributed between the branch and the cycle;
when $N$ is odd, the black nodes of the cycle are one more than those in the branch.
Finally, here is a picture of the dual graph of  $\mathrm{Dl}S=[s_{k_1}s_{k_2}...s_{k_N}r_l]$
in the case \underline{$N$ even} and $l\geq 2$. The $N=odd$ case has been shown in \cite{fp15} with a slightly different notation. In the case $l=1$ the root of the tree is the first cycle black node $C_0$.

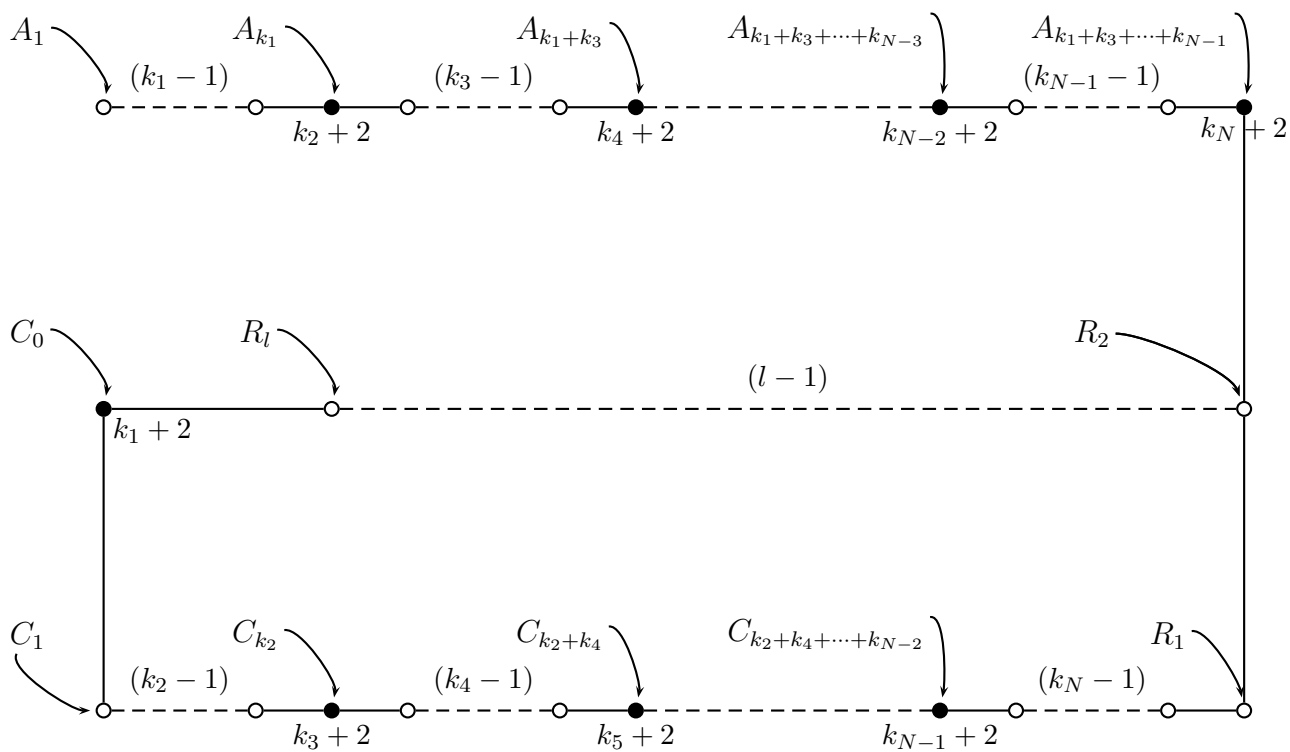
\begin{pspicture}(2,0)(20,12)

\cnode(2,10){3 pt}{A}
\rput(1,11){\rnode{a}{$A_1$}}

\cnode(4,10){3 pt}{B}

\cnode[fillstyle=solid,fillcolor=black](5,10){3 pt}{C}
\uput[d](5,10){\small{$k_2+2$}}
\rput(4,11){\rnode{c}{$A_{k_1}$}}

\cnode(6,10){3 pt}{D}

\cnode(8,10){3 pt}{E}

\cnode[fillstyle=solid,fillcolor=black](9,10){3 pt}{F}
\uput[d](9,10){\small{$k_4+2$}}
\rput(8,11){\rnode{f}{$A_{k_1+k_3}$}}

\cnode[fillstyle=solid,fillcolor=black](13,10){3 pt}{G}
\uput[d](13,10){\small{$k_{N-2}+2$}}
\rput(11.5,11){\rnode{g}{$A_{k_1+k_3+\cdots+k_{N-3}}$}}

\cnode(14,10){3 pt}{H}

\cnode(16,10){3 pt}{I}

\cnode[fillstyle=solid,fillcolor=black](17,10){3 pt}{J}
\uput[dr](16.3,10){\small{$k_N +2$}}
\rput(15.5,11){\rnode{j}{$A_{k_1+k_3+\cdots+k_{N-1}}$}}

\ncline[linestyle=dashed]{A}{B}
\Aput{\small{$(k_1-1)$}}
\ncline{B}{C}
\ncline{C}{D}
\ncline[linestyle=dashed]{D}{E}
\Aput{\small{$(k_3-1)$}}
\ncline{E}{F}
\ncline[linestyle=dashed]{F}{G}
\ncline{G}{H}
\ncline[linestyle=dashed]{H}{I}
\Aput{\small{$(k_{N-1}-1)$}}
\ncline{I}{J}

\nccurve[ncurv=.4,angleB=10,angleA=80,nodesep=2pt]{<-}{A}{a}
\nccurve[ncurv=.4,angleB=10,angleA=80,nodesep=2pt]{<-}{C}{c}
\nccurve[ncurv=.4,angleB=10,angleA=80,nodesep=2pt]{<-}{F}{f}
\nccurve[ncurv=.4,angleB=10,angleA=80,nodesep=2pt]{<-}{G}{g}
\nccurve[ncurv=.4,angleB=10,angleA=80,nodesep=2pt]{<-}{J}{j}

\cnode[fillstyle=solid,fillcolor=black](2,6){3 pt}{A''}
\uput[dr](2,6){\small{$k_1+2$}}
\rput(1,7){\rnode{a''}{$C_0$}}

\cnode(5,6){3 pt}{B''}
\rput(4,7){\rnode{b''}{$R_l$}}

\cnode(17,6){3 pt}{C''}
\rput(15,7){\rnode{d''}{$R_2$}}
\nccurve[ncurv=.4,angleB=0,angleA=110,nodesep=2pt]{<-}{C''}{d''}

\ncline{A''}{B''}
\ncline[linestyle=dashed]{B''}{C''}
\Aput{\small{$(l-1)$}}

\nccurve[ncurv=.4,angleB=10,angleA=80,nodesep=2pt]{<-}{A''}{a''}
\nccurve[ncurv=.4,angleB=10,angleA=80,nodesep=2pt]{<-}{B''}{b''}
\nccurve[ncurv=.4,angleB=0,angleA=110,nodesep=2pt]{<-}{C''}{d''}

\cnode(2,2){3 pt}{A'}
\rput(1,3){\rnode{a'}{$C_1$}}

\cnode(4,2){3 pt}{B'}

\cnode[fillstyle=solid,fillcolor=black](5,2){3 pt}{C'}
\uput[d](5,2){\small{$k_3+2$}}
\rput(4,3){\rnode{c'}{$C_{k_2}$}}

\cnode(6,2){3 pt}{D'}

\cnode(8,2){3 pt}{E'}

\cnode[fillstyle=solid,fillcolor=black](9,2){3 pt}{F'}
\uput[d](9,2){\small{$k_5+2$}}
\rput(8,3){\rnode{f'}{$C_{k_2+k_4}$}}

\cnode[fillstyle=solid,fillcolor=black](13,2){3 pt}{G'}
\uput[d](13,2){\small{$k_{N-1}+2$}}
\rput(11.5,3){\rnode{g'}{$C_{k_2+k_4+\cdots +k_{N-2}}$}}

\cnode(14,2){3 pt}{H'}

\cnode(16,2){3 pt}{I'}


\cnode(17,2){3 pt}{J'}
\rput(16,3){\rnode{j'}{$R_1$}}


\ncline[linestyle=dashed]{A'}{B'}
\Aput{\small{$(k_2-1)$}}
\ncline{B'}{C'}
\ncline{C'}{D'}
\ncline[linestyle=dashed]{D'}{E'}
\Aput{\small{$(k_4-1)$}}
\ncline{E'}{F'}
\ncline[linestyle=dashed]{F'}{G'}
\ncline{G'}{H'}
\ncline[linestyle=dashed]{H'}{I'}
\Aput{\small{$(k_N-1)$}}
\ncline{I'}{J'}

\ncline{A'}{A''}
\ncline{J}{C''}
\ncline{J'}{C''}

\nccurve[ncurv=.4,angleB=240,angleA=180,nodesep=2pt]{<-}{A'}{a'}
\nccurve[ncurv=.4,angleB=10,angleA=80,nodesep=2pt]{<-}{C'}{c'}
\nccurve[ncurv=.4,angleB=10,angleA=80,nodesep=2pt]{<-}{F'}{f'}
\nccurve[ncurv=.4,angleB=10,angleA=80,nodesep=2pt]{<-}{G'}{g'}
\nccurve[ncurv=.4,angleB=10,angleA=100,nodesep=2pt]{<-}{J'}{j'}

\end{pspicture}

\bigskip

The following multilinear forms have been introduced in \cite[3.1]{fp15} and are tailored to satisfy \cite[Corollary 2.7]{fp15}. They play an important role in the description of several logarithmic invariants of the pair $(S,D)$.

\begin{Definition} \rm Let $X_1, X_2,...,X_n$ denote a set of variables
and define polynomials $f$ in $n$ variables inductively, by
$$ f_{-1}:=0,  \quad f_0:=1,  \quad f(X_1):= X_1,   \quad f(X_1,X_2):=X_1X_2 +1 $$
\be f(X_1,...,X_n):=X_n f(X_1,...,X_{n-1}) + f(X_1,...,X_{n-2}).
\label{f} \ee
 We also consider multilinear forms $\cp$, inductively defined from $f$ as follows
 $$  \mathcal{P}(X_1) = X_1, \quad  \mathcal{P}(X_1,X_2)=X_1X_2+X_1 $$
\be \cp(X_1,...,X_n) := X_n  f(X_1,...,X_{n-1}) + \cp(X_1,...,X_{n-1})
\label{p}\ee
 \hfill     $\bigtriangleup$  \end{Definition}

\begin{Remark} \rm
It is easy to check that
$$ \cp(X_1,...,X_n) +1 =  f(X_1,...,X_n) +  f(X_1,...,X_{n-1})  $$
  \hfill $\bigtriangleup$
\end{Remark}
\medskip

Our first relation deals with the coefficients of the following $2\times 2$~matrix introduced by Dloussky. It encodes the sequence of non-generic $(b-l)$ blow ups describing $S$ as a surface with a global spherical shell, see \cite[p.84, p.91]{dl16}.

$$ 
\left( \begin{array}{cc}
p & q  \\
r & s   \end{array} \right) :=
\left( \begin{array}{cc}
0 & 1  \\
1 & k_1   \end{array} \right)  \cdots
\left( \begin{array}{cc}
0 & 1  \\
1 & k_{N}   \end{array} \right)
$$

It is in fact easy to check by induction that

\begin{Remark} \rm The coefficients of Dloussky matrix are given by the multilinear polynomials $f$:
$$\begin{array}{cc}
s= & f(k_1,\dots,k_{N}) \\
r= & f(k_1,\dots,k_{N-1}) \\
q= & f(k_2,\dots,k_{N}) \\
p= & f(k_2,\dots,k_{N-1})
\end{array}$$
In particular,
\be r+s = \cp(k_1,...,k_{N}) +1 \quad \textrm{ and } \quad p+q = \cp(k_2,...,k_{N}) +1 \ee \label{ks}   
\noindent and from \cite[5.3]{fp15} we recover the result of \cite[4.24]{dl16} that $r+s$ equals the topological invariant $k(S)$ 
which equals the opposite of the determinant of the intersection matrix of the irreducible components $\{A_1,\dots,A_{k_{1}+\dots k_{N-1}}\}$ of the unique tree (see the picture above) \cite[thm. 3.20]{dl11}.

  \hfill $\bigtriangleup$
\end{Remark}
\medskip

\section{First Chern class and birational germs}

Except for Enoki surfaces where the cycle $D$ is trivial in homology,  the rational curves of a Kato surface form a basis of  $H^2(S,\q)$ and in \cite{fp15} we solved the problem of computing the rational coefficients of the anti-canonical class,  namely the first Chern class $c_1 \in H^2(S,\z)$ for every intermediate Kato surface, with arbitrary number of trees. 
The aim of this section is to discuss some connections between our solution and new birational structures introduced by Dloussky for simple intermediate Kato surfaces. For this reason we will recall our results in this particular case in which $\Gamma_D$ has a unique tree.

Let $b:=b_2(S)$ and $D=\sum_{i=0}^b D_i$ be the maximal reduced divisor of an intermediate Kato surface, from the previous discussion there exist unique rational coefficients $d_i\in\q$ which we call \it multiplicities \rm 
such that the following equation holds in cohomology 
$$c_1 = \sum_{i=1}^{b_2(S)}d_jD_j  \in H^2(S,\q)$$

This equation is completely equivalent to the linear system \cite[(1) p.1532]{do99} via the adjunction formula applied to every irreducible component $D_i(\cong\mathbb{CP}_1)$ of the maximal curve $D$: \be 2=(c_1 -D_i)D_i. \ee

Setting $\t $ to be the multiplicity $a_1$ of the tip $A_1$ of the unique tree -- see the picture of the previous section -- we proved that $\t$ is always strictly positive and smallest among the $d_j$'s which can all be expressed in terms of the tip multiplicity $\t$ as follows: by \it White Lemma \rm \cite[2.1]{fp15} the multiplicities grow up linearly along chains of $(-2)$-curves (white nodes) and for the chain $R_2\cdots R_l$ we showed the perhaps surprising result that the slope is always $-1$. Namely, for every $S$:  
\be r_l=\t+2;\cdots ; r_2=\t + l.   \ee \label{regular}

Again surprisingly, the multiplicity $c_0$ of the first cycle black node $C_0$ is always $\t+1$ -- for every $S$ -- and in general the growth of the $d_i$ is  piecewise linear with slope $g_j$ changing at every black node $j=1,\dots,N+1$  as described in  \it Black Lemma \rm   \cite[2.1]{fp15};  using the duality between black nodes in the cycle and white nodes in the branch (and viceversa) we showed that 
\be  g_{j_+1}-g_{j-1}=k_{j+1}g_j  \ee   
\cite[Corollary 2.7]{fp15} -- in the notations of the present paper. This formula shows that the multilinear polynomials $f(X_1,\dots,X_n)$ are suited to express all multiplicities and setting as in \cite[4.24]{dl16}
  $$ 
\left( \begin{array}{cc}
p_j & q_j  \\
r_j & s_j   \end{array} \right) :=
\left( \begin{array}{cc}
0 & 1  \\
1 & k_1   \end{array} \right)  \cdots
\left( \begin{array}{cc}
0 & 1  \\
1 & k_j   \end{array} \right)
$$
we showed that the black nodes multiplicities are always equal to $1+g_j$; furthermore, for every $j=1,\dots , N$
\be g_{j+1}=s_j\t-q_j \ee           \label{slopes}
\noindent and this allowed us to express the root multiplicity 
\be  r_2=g_{N+1} +g_N +1 = s\t -q +r\t -p +1=(r+s)\t-(p+q)+1. \ee \label{root}
\noindent on the other hand we also have $r_2=\t+l$ by (\ref{root}) and therefore we conclude that  \cite[3.4]{fp15}   
\be   \t=\frac{p+q-1+l}{r+s-1}   \ee   \label{index}

Recall now the notion of index:

\begin{Definition} \rm   The \it index \rm of a Kato surface $S$ is the smallest natural number $d$ such that the cohomology class $d \cdot c_1(S)$ can be expressed as a linear combination of the rational curves $D_i$ with \it integer \rm coefficients and is therefore represented by an effective divisor. In other words $d$ is the smallest (automatically positive) integer such that $H^0(S, \mathcal K ^{-d}\otimes F)\cong \c$ for some $F\in $ Pic$_0(S)$ where $\mathcal K$ denotes the canonical line bundle of $S$ and Pic$_0(S)\cong\c^*$ is the group of flat holomorphic line bundles.
\hfill $\bigtriangleup$
\end{Definition}

We then have that $\mathrm{index}(S)$ is the least common multiple of all denominators of the rational coefficients of the first Chern class $c_1(S)=\sum_{j=1}^b d_i D_i$. Because we have just expressed all multiplicities as linear functions of $\t$ with integer coefficients we can recover from the above expression of $\t$ the result of \cite[4.25]{dl16} that index$(S)=g.c.d.\{ p+q-1+l , r+s-1\}$.

\bigskip

An important tool in the study of Kato surfaces was introduced by Dloussky \cite{dl84} who showed that the complex structure of a Kato surface $S$ is completely determined by the contracting germ of the map $\pi\circ\psi$ at the origin $0\in B$, where $\pi$ is the blowing down map. Polynomial normal forms for this germ were given by Favre \cite{fa00}; more recently Dloussky \cite{dl16} introduced new germs which for a simple intermediate Kato surface $S$ with Dloussky sequence DlS$=[s_{k_1}s_{k_2}...s_{k_N}r_l]$ take the following form, we will call them \it birational germs \rm:
$$G(z_1,z_2)=(z_1^{p+rl}z_2^{q+sl}+\sum_{1=0^{l-1}}a_1(z_1^rz_2^s)^{i+1}+a_{l+K}(z_1^rz_2^s)^{l+K+1},z_1^rz_2^s)$$
if the first blow up is not generic. Otherwise, when the first blow up is generic 

$$G(z_1,z_2)=((z_1z_2^l +\sum_{i=0}^{l-1}a_iz_2^{i+1} + a_{l+K} z_2^{l+K+1} )^p z_2^q ,
                          (z_1z_2^l +\sum_{i=0}^{l-1}a_iz_2^{i+1} + a_{l+K} z_2^{l+K+1} )^r z_2^s)$$ 
where $(z_1,z_2)$ are coordinates around $0\in\c^2$ and $p,q,r,s,l$ are the integers of the previous section which encode the sequence of blowups and the length of the regular sequence.

The aim of this section is to highlight some (new) interactions between these germs and the first Chern class $c_1(S)\in H^2(S,\z)$.

We start by recalling that the integer $K$ in the birational germ is defined to be 
\be K=max\{0, [\frac{l-d}{r+s-1}]\} \ee
\noindent  \cite[4.27, 4.23]{dl16} for 
\be d:=(r+s)-(p+q) \ee \label{d}
\noindent  while the rational numbers 
\be u:=\frac{p+s+rl-1-(1)^N}{r+s-1} \quad\mathrm{ and }\quad v:= \frac{q+r+sl-1+(1)^N}{r+s-1} \ee\label{u}
\noindent  introduced in \cite[4.22]{dl16} are the vanishing order of a twisted anti-canonical section along the axes $z_2=0$ and $z_1=0$.

Before stating our result let us observe that in an intermediate Kato surface each rational curve meets exactly two other rational curves except for the tip of the tree $A_1$ which meets only one and for the root $R_2$ which meets three, more precisely, when $N$ is even: $R_1$, $A_{k_1+\cdots + k_{N-1}}$ (black node) and $C_{k_2+\cdots +k_N}$ (white node).

\begin{Proposition}\label{t-1}
The rational number  $\frac{l-d}{r+s-1}$ equals the tip multiplicity minus $1$: $\frac{l-d}{r+s-1}=\t-1$.

The rational numbers $u$ and $v$ coincide with the multiplicities of the curves $R_1$ and of the black node meeting the root $R_2$, respectively. 

In particular, $\frac{l-d}{r+s-1}$ is an integer if and only if $S$ has index $1$ if and only if $u$ and $v$ are both integers \cite[4.22]{dl16}. 
\end{Proposition}
{\it Proof.} Of course, $\frac{l-d}{r+s-1}= \frac{p-r-s+p+q}{r+s-1}=\frac{l+p+q-1}{r+s-1}-1=\t-1.$ 

\noindent In a similar vein, it is easy to check that $u=g_N +1=r\t-p+1$ while $v=g_{N+1} +1=s\t-q+1$  by (\ref{slopes}) and therefore are both integers when $\mathrm{index}(S)=1$. Viceversa, $u+v=g_N + g_{N+1} +2=r_2+1$ by the result of (\ref{root})
which is also equal to $\t +l +1$ and is therefore an integer if and only if $\t$ is an integer -- i.e. $\mathrm{index}(S)=1$. 
\hfill $\bigtriangleup$
\bigskip

\section{Deformations}

In this section we study the deformation theory of intermediate Kato surfaces, find the dimension of the Kuranishi family by computing the cohomology of a relative tangent bundle and discuss the relation with the number of parameters in the birational forms of the preceding section. 

While birational forms exist only for intermediate surfaces with a single tree the cohomology computation holds in general for any Kato surface and to this end we recall that a general Kato surface $S$ of intermediate type has Dloussky sequence of the form DL$S= [$Dl$S_1 \ldots $ Dl$S_F$] where each Dl$S_f = [s_{k_{1}f} \ldots s_{k_{N}f}r_{lf}]$ is simple and $F$ denotes the number of trees in the associated dual graph $\Gamma_D$ which is constructed from DL$S$ using the same rules as before. 

By Riemann-Roch the Euler characteristic of the tangent bundle $\Theta_S$ turns out to be always $2b$ -- twice the second Betti number of $S$. From a vanishing theorem of Nakamura it then follows that $h^1(\Theta_S)=2b + h^0(\Theta_S)$ is the dimension of the deformations of $S$ as a complex surface -- of course in class VII -- with generic member of the family a blown-up Hopf surface. Because we are interested in the deformations of $S$ among Kato surfaces we consider deformations of the pair $(S,D)$ as in \cite[2.1]{fp10}; any element $S_t$ in this deformation family will have a fixed configuration of rational curves with dual graph equal to $\Gamma_D$ and is therefore automatically minimal with $b$ rational curves; by \cite{dot03}  Kato surfaces are characterized, among surfaces in class VII$_0$, by the property of having $b=b_2(S)$ rational curves; this shows that $S_t$ is a Kato surface with the same Dloussky sequence of $S$.

The deformations of the pair $(S,D)$ are governed by the cohomology of the relative tangent bundle $\Theta_S(-\log D)$ of vector fields which are tangent to $D$, along $D$. We will show that its cohomology only depends on Dl$S$ -- more precisely, on the total length of its regular sequences -- except for the case in which $S$ has a holomorphic vector field, automatically tangent to $D$, along $D$; this can only occur when index$(S)=1$, a necessary but not sufficient condition \cite{do99}.
\medskip

Let $S$ be an intermediate Kato surface, its maximal curve $D$ has $b$ irreducible components $D_i, 1\le i\le b$ 
 and let Dl$S$ = [Dl$S_1 \ldots $ Dl$S_F$] be its Dloussky sequence with $F$ the number of trees of the dual graph $\Gamma_D$.  
In the following proposition we directly compute the dimension of the tangent space $H^1(S,\Theta_S(-\log D))$ of the Kuranishi family of deformations of the pair $(S,D)$. We set
$\epsilon = \epsilon_S:= \dim H^0(S,\Theta_S)=H^0(S,(-\log D)$ and recall that $\epsilon\leq 1$ where equality implies 
index$(S)=1$.

\begin{Proposition}\label{}  \cite{dl14}
For any Kato surface $S$, $\dim H^1(S,\Theta_S(-\log D))= l + \epsilon$ where $l=l1+\dots +lF$ is the total number of irreducible components in the regular sequences of Dl$S$.
\end{Proposition}

{\em Proof.}
We consider the sheaf exact sequence
\begin{equation}\label{}
 0 \ra \Theta_S(-\log D) \ra \Theta_S  \ra N_D \ra 0   
\end{equation}
where $N_D=\oplus_{i=1}^b N_i$ is the normal bundle of the normal-crossing divisor $D$ and each $N_i = N_{D_i/S}\cong \mathcal O _{\mathbb{P}_1}(D_i^2)$ is the normal bundle of $D_i$ in $S$. Suppose first that all $D_i$ are smooth. Then, since $D^2_i$ is negative we have $h^0(N_i)=0$ for all $i$ and
we recover that $h^0(-\log\Theta)=h^0(\Theta)$. Furthermore, by a result of Nakamura \cite[\S 3]{na90}, see also \cite{fp10}, we have $h^2(\Theta_S)=h^2(\Theta_S(-\log D))=0$ so that $\chi(\Theta_S(-\log D))=2b -\sum_{i=0}^b h^1(N_i)$. Now, $h^1(N_i)=1$ for each white node in $\Gamma_D$ while for a black node $D_j^2=-(k_j+2)$ 
we will have $h^1(N_j)=k+1$. If $N$ denotes as usual the number of balck nodes we have that the number of white nodes is $b-N$ while $\sum_{j=0}^N k_j = b-l$ and we conclude that $h^1(N_D)=b-N+b-l+N=2b-l$ and the result follows in this smooth case.

It remains to discuss the exceptional case when $S$ has a singular rational curve which happens exactly for Dl$S=[s_{b-1}r_1]$ and therefore the cycle has a unique component, singular with a double point. We consider in this case the unique unramified double cover $u:\tilde S \to S$; it will have Dl$\tilde S= [s_{b-1}r_1s_{b-1}r_1]$ and furthermore $\Theta_{\tilde{S}}(-\log \tilde{D}))\cong u^*\Theta_S(-\log D))$ therefore by Riemann-Roch it follows that 
$$\chi(\tilde{S},\Theta_{\tilde{S}}(-\log \tilde{D})) = 2\chi(S,\Theta_S(-\log D))=2l=2$$ 
from the first part of the proof.  We already know that the second cohomology always vanishes, furthermore $\epsilon(\tilde S)=\epsilon(S)$ and we can easily conclude that $h^1(\Theta(-\log D)))=1+\epsilon$ also in this exceptional case $l=1$; notice that index$([s_{b-1}r_1])=1$ implies $b=2$.
\hfill  $\Box$
\medskip

It is perhaps worthwhile to point out that the above result holds true for any Kato surface, with the same proof. It says that hyperbolic and half-Inoue surfaces, corresponding to the $(l=0)-$case, are logarithmically rigid while Enoki surfaces, the opposite case $b=l$, have a $b$-dimensional Kuranishi family of log-deformations.

\medskip
 
\begin{Remark}  \rm
The above result is related to the number of parameters in the birational germ of intermediate Kato surfaces with one branch. It is shown by Dloussky that the parameters $a_0,\dots , a_{l-1}$ are always effective while the extra parameter $a_{l+K}$ is effective only for surfaces with a holomorphic vector field -- i.e. $\epsilon=1$ -- and can otherwise be taken to vanish. These considerations extends to all the other Kato surfaces which also always have birational structures \cite{dl16}.

The fact that the dimension of the logarithmic moduli space equals the total length of the regular sequences is geometrically explained as follows.  In the correspondence between the blowing-up sequence and the Dloussky sequence  each entry of a singular sequence corresponds to a blowing up with center one of the nodes of the previously produced exceptional curves and hence has no moduli, while for an entry of a regular sequence the blowing up occurs at a general point and hence each contributes to one dimensional moduli.
\hfill  $\bigtriangleup$    \end{Remark}

We can now compute the dimension of the tangent space $H^1(S, \Theta_S(-D))$ of the Kuranishi family of deformations of $S$ vanishing along the maximal curve $D$. 

\begin{Proposition} Let $S$ be a Kato surface, then $\dim H^1(S, \Theta_S(-D)) = b+l-\eta$. Where $b$ is the second Betti number of $S$, $l$ is the total number of components in the regular sequences of $S$ while $\eta=1$ if and only if there is a holomorphic vector field and a tip of multiplicity $\t=1$. In particular, for simple intermediate Kato surfaces we have that $\eta=1$ implies that there is a vector field, the length of the regular sequence satisfies $l=r+s-(p+q)$ and this holds if and only if in the birational germ the coefficient $a_{l+K}=a_l$ is effective.
\end{Proposition}
\it Proof. \rm 
Consider the exact sequence of sheaves in $S$
$$0\to \Theta_S(-D) \to \Theta_S(-\log D) \to \Theta_D \to 0$$ where $\Theta_D$ is the tangent sheaf of the maximal curve $D$. Since $D$ is the union of $b$ irreducible components $D_i$ meeting transversally and which we can assume to be smooth rational curves, we see that $\Theta_D=\oplus_{i=1}^b T_i$ with each 
$T_i\cong  \Theta_{{\mathbb P}_{1}}(-a_{i})$ where $a_i$ denotes the number of irreducible components intersecting $D_i$.  

Therefore, $T_i\cong \mathcal{O}_{\mathbb{CP}_1}$ most of the times while $T_i\cong \mathcal{O}_{\mathbb{CP}_1}(1)$ if $D_i$ is a tip and the only other possibility is that $D_i$ is a root in which case $T_i\cong \mathcal{O}_{\mathbb{CP}_1}(-1)$. It then follows that $h^1(\Theta _D)=0$ and since the number of tips equals the number of roots we also have that  $h^0(\Theta _D)=b$, always.

 From the previous Proposition we easily get that the Kuranishi
family is unobstructed with $\quad \dim H^1(S, \Theta_S(-D)) = b+l \quad$ 
whenever the following map is isomorphic:
$H^0(S, \Theta_S(-D)) \to H^0(S,\Theta_S(-\log D))$.

The only other possibility is that $S$ has a holomorphic vector field, automatically tangent to $D$ along $D$ and not vanishing identically on $D$ in which case we have to subtract $\eta=1$ dimensions and notice that index$(S)=1$ so that $c_1(S)$ is represented by the effective divisor $\sum_{i=0}^b d_i D_i$ whose coefficients $d_i$ we computed in the previous section. Finally, it is shown in \cite{do99} that a holomorphic vector field $\theta$ vanishes exactly on the divisor $D_{\theta} :=  (\sum_{i=0}^b d_i D_i) -D = \sum_{i=0}^b (d_i -1) D_i$. Therefore $\theta$ is not identically zero on the maximal curve $D$ if and only if $1$ is the minimal multiplicity of $c_1(S)$. But we have already shown that the absolute minimum for the multiplicities of the first Chern class of an intermediate Kato surface occurs at one of its tips, so that $\t=1$; therefore $l=d$ by (\ref{index}) and $K=0$.
\hfill $\Box$

\begin{Remark}  \rm
Every $S$ with a holomorphic vector field and Dl$S=[s_k r_k]$ gives an example of the above situation $\eta=1$.
\hfill        $\bigtriangleup$        \end{Remark}

 \bigskip
 
\noindent \it Twisting line bundles. \rm
 
 For any Kato surface $S$ we have Pic$_0(S)\cong\c^*$ because holomorphic flat line bundles are parametrized by representations of $\pi_1(S)\cong\z$; after fixing a generator of $\pi_1(S)$ we write elements in Pic$_0(S)\cong\c^*$ as $L^\gamma$ with $\gamma=1$ corresponding to the trivial complex line bundle $\mathcal O _S$.
Using this notation is was proved in \cite{do99} that a fixed intermediate Kato surface $S$ with canonical bundle $\mathcal K$ has index$(S)=1$ if and only if there is a unique $\alpha\in\c^*$ such that $H^0(S,\mathcal K^{-1} \otimes L^\alpha)\neq 0$ if and only if there is a unique $\beta\in\c^*$ such that $H^0(S, \Theta \otimes L^\beta)\neq 0$.
By \cite[4.34]{dl16} the relation between these exponents is given by $\beta ^{-1} = k(S)\alpha$ where $k(S)=r+s$ is the topological invariant of (\ref{ks}).  Furthermore, by \cite[4.22, 4.35]{dl16} the twisting coefficient $\alpha$ only depends on the first parameter $a_0$ of the birational germ in such a way that
 \be \alpha=(-1)^N a_0^u \ee
 where $u$ is the natural number of (\ref{u}) and $N$ is the number of black nodes. In particular, as $a_0$ varies in $\c^*$ we see that given any $\gamma\in\c^*$ there is an intermediate Kato surface $S$ with $H^0(S,\mathcal K^{-1} \otimes L^\gamma)\neq 0$. 
 
We now want to use the birational germ to distinguish among three interesting subclasses of intermediate Kato surfaces of index $1$; we set $H_{vf}$ to be the surfaces with vector fields and $H_{ac}$ those with anti-canonical section, while $H_{bh}$ denotes those surfaces which can possibly admit a bi-Hermitian metric. The following result was proved in \cite[5.11]{fp15}
 
 \begin{Proposition} The three subsets of intermediate Kato surfaces of index $1$ given by $H_{vf}$, $H_{ac}$ and $H_{bh}$ are mutually disjoint.
\end{Proposition}
 {\it Proof.} We look at the coefficient $a_0$ in the birational germ which is an effective parameter of deformations. Clearly $|a_0(S)|=1$ if $S\in H_{ac}$ while for $S\in H_{vf}$ we will have $H^0(S,K^{-1}\otimes L^\alpha)\neq 0$ with $\alpha=1/(r+s)$ so that $|a_0(S)|<1$ in this case because both $r+s$ and $u$ are certainly $>1$. Finally, if $S\in H_{bh}$ we have $H^0(S,K^{-1}\otimes L^\alpha)\neq 0$ with deg$(L^\alpha)<0$ by \cite{ap01}. According to \cite[4.21]{dl14} the Gauduchon degree of the holomorphic line bundle $L^\gamma$, whose sign does not depend on the choice of a particular Gauduchon metric because $b_1(S)=1$, has the opposite sign of $\log |\gamma|$. Therefore we have that $|a_0(S)| >1$ when $S\in H_{bh}$.
 \hfill $\Box$.  \medskip

We conclude with the simple observation that, to the contrary, the bi-Hermitian parabolic Inoue surfaces produced in \cite{po97}, \cite{fp10} 
admit both a holomorphic vector field and anti-canonical section.


\newcommand{\bysame}{\leavevmode\hbox to3em{\hrulefill}\,}

\bigskip
\noindent
M. Pontecorvo --  Dipartimento di Matematica e Fisica, Roma Tre University.

\end{document}